\documentclass{article}
\usepackage{graphicx}
\usepackage{amsmath}
\usepackage{amsthm}
\usepackage{caption}
\providecommand{\keywords}[1]{\textbf{Keywords} #1}
\newtheorem{remark}{Remark}
\newtheorem{theorem}{Theorem}

\begin{document}

\title{A priori and a posteriori error estimation for finite element approximation of advection-diffusion-reaction equation with spatially variable coefficients}
\author{Prof. B.V. Rathish Kumar \thanks{drbvrk11@gmail.com} \\  Manisha Chowdhury \thanks{chowdhurymanisha8@gmail.com} }
      
\date{Indian Institute of Technology Kanpur \\ Kanpur, Uttar Pradesh, India}

\maketitle

\begin{abstract}
This paper presents a study of finite element error estimation of advection-diffusion-reaction equation with spatially variable coefficients. We have derived a priori and a posteriori errors in both energy and $L_2$ norm. We have used residual-based a posteriori error estimator. Numerical results are also presented to verify our theoretical approach.
\end{abstract}


\section{Introduction}
\hspace{3mm} Dispersion, an important phenomena in fluid dynamics, describes the spreading of mass from highly concentrated areas to less concentrated areas. Solute dispersion in a straight tube has broad applications in the fields of bio-medical engineering, chemical engineering, environmental science etc. Some specific real life examples are: transportation of drug or toxins through body fluid, use of drug eluting stents into blood vessels, transportation of pollutant in environment etc. \\

All these physical phenomena  can be adequately described mathematically by the advection-diffusion-raction (ADR) equation. Due to complexity in many of these problems, the ADR equation is solved numerically. Numerical methods provides an approximate solution to the physical problem. Therefore the convergence analysis of the approximate solution to the exact solution is an important topic to study. There are several studies on error estimation with constant coefficients (\cite{RefB}-\cite{RefE},\cite{RefG}-\cite{RefJ}), but physical problems involve some difficulties such as variable coefficients, non-uniform flow etc. Here we are dealing with ADR equation with spatially variable diffusion and velocity coefficients for steady, incompressible flow. This kind of ADR equation plays an important role in studying pollutant transportation in the field of environmental science. There are very few works in finding analytical solution of ADR equation with spatially variable coefficients \cite{RefF}, but they studied a particular form. As per our knowledge, no study has been done on convergence analysis of a numerical scheme, finding approximate solution of the ADR equation with spatially variable coefficients.  Here we have found out the bounds of errors obtained in a priori and a posteriori estimates under finite element method (FEM). FEM is a well-known numerical scheme for solving partial differential equations of boundary value problems in the areas of structural analysis, heat and mass transfer, fluid flow etc. Both the estimation give order of convergence of the error, where a priori shows stability of the method, the study of a posteriori estimate is useful for adaptivity and control of error solution. We have here discussed only the theoretical aspect of a posteriori error  estimation.\\  

This paper is organised as follows: second section consists of various of things starting from defining the problem, its strong and weak form, existence and uniqueness of the solution of variational form and Galerkin finite element formulation. Third section establishes stability estimates and in the next four sections, a priori and a posteriori error estimations in both energy and $L_2$ norm are presented. At the last section we have introduced the figures and results in tabular form obtained from numerical computations to verify the order of accuracy in error estimates.
\section{Statement of the problem}

\hspace{2mm}The two-dimensional solute transport equation with spatially variable dispersion coefficients on a bounded domain $\Omega \subset R^2 $ with homogeneous Dirichlet boundary condition for a steady flow can be stated as follows:\\
Find c : $ \Omega \rightarrow $ R such that\\
\begin{equation}
     \frac{\partial}{\partial x}[D_x \frac{\partial c}{\partial x}- u c] + \frac{\partial}{\partial y}[D_y \frac{\partial c}{\partial y} - v c] - \mu c + q = 0 \hspace{2mm} in \hspace{1mm} \Omega
\end{equation}
\hspace{8.2cm} c = 0 on $\partial\Omega $\vspace{2mm}\\ where c is the solute concentration of the dispersing mass, $(D_x,D_y)$ and $(u,v)$ are the diffusion coefficients and the velocity components along x-axis and y-axis respectively at the position $(x,y)$, $\mu$ is the reaction coefficient and q represents source term.
\label{sec:1}
\subsection{Incompressibility condition}
For incompressible fluid, the divergence of advection velocity is zero that is,\\

 $\nabla \cdot \overline{u}$ = 0 on $\Omega$, where $\overline{u} = (u,v)$
 
 \subsection{Strong form}
 Notations $\widetilde{\bigtriangledown}$ := $(D_x \frac{\partial}{\partial x},D_y \frac{\partial}{\partial y})$ \vspace{2mm}\\
Applying the incompressibility condition on (1), the strong form of the problem is as follows:\\ Find c $\in H^2(\Omega)$ such that\\
\begin{equation}
\mathcal{L}c = -\bigtriangledown \cdot \widetilde{\bigtriangledown}c + \overline{u} \cdot \bigtriangledown c + \mu c = q \hspace{2mm} in \hspace{1mm} \Omega
    \end{equation} 
\hspace{7cm} c = 0 on $\partial\Omega $\vspace{2mm}\\
Let us take source term q $\in  L_2(\Omega)$.\\
Now we introduce the spaces mentioned above along with their respective norms as follows:\vspace{3mm}\\
$L_2(\Omega) = \{v \in \Omega \mid \int_{\Omega} v^2 d\Omega < \infty \}$ \vspace{3mm}\\
$H^2(\Omega)= \{ v \in \Omega \mid D^\alpha v \in L^2(\Omega) ; \forall \alpha \hspace{1 mm} $such that$ \hspace{1mm} 0\leq \mid \alpha \mid \leq 2 \}$ \vspace{3mm}\\
where $D^\alpha v = \frac{\partial^{\mid \alpha \mid} v}{\partial x ^{\alpha_1}\partial y^{\alpha_2}} \hspace{2mm} $for each$ \hspace{2mm} \alpha = (\alpha_1,\alpha_2), where \hspace{1mm} \alpha_i \geq 0, \hspace{1mm}\alpha_i=integer, \hspace{1mm} for \hspace{1mm} i=1,2 \hspace{1mm} and  \mid \alpha \mid = \alpha_1 + \alpha_2$ \vspace{3mm}\\
let $\| \cdot \|$ and $\| \cdot \|_2$ denote $L_2$ and $H^2$ norms respectively over $\Omega$. They are defined as follows:\\
$\|v\|=(\int_{\Omega} v^2 )^{1/2}$ and \\
$\|v\|_2 = (\|v\|^2 +\| \frac{\partial v}{\partial x} \|^2 + \|\frac{\partial v}{\partial y} \|^2 + \|\frac{\partial^2 v}{\partial x^2} \|^2 + \|\frac{\partial^2 v}{\partial x \partial y} \|^2 + \|\frac{\partial^2 v}{\partial y^2} \|^2)^{1/2}$\\

\subsection{Weak form}
In order to write weak form of equation (2) with homogeneous Dirichlet boundary condition let us first introduce the space, \vspace{3mm}\\
$ H^1 (\Omega)$ = $\{ v \in \Omega \mid D^\alpha v \in L^2(\Omega) ; \forall \alpha \hspace{1 mm} $such that$ \hspace{1mm} 0\leq \mid \alpha \mid \leq 1 \hspace{1mm} \}$ \vspace{3mm}\\
let $\| \cdot\|_1$ denote $H^1$ norm over $\Omega$ and it is defined by, \vspace{2mm}\\
$\|v\|_1 = (\|v\|^2 +\| \frac{\partial v}{\partial x} \|^2 + \|\frac{\partial v}{\partial y} \|^2)^{1/2}$ \vspace{2mm}\\
let V:= $H_0^1(\Omega)=\{v \in H^1(\Omega) | \hspace{1mm} v=0 \hspace{1mm} on \hspace{1mm} \Omega \}$ \vspace{2mm}\\
The weak form be: Find $c\in V$ such that 
\begin{equation}
  a(c,d)=l(d) \quad \forall d \in V  
\end{equation}
where a(c,d)=$\int_{\Omega}D_x\frac{\partial c}{\partial x}\frac{\partial d}{\partial x} + \int_{\Omega}D_y\frac{\partial c}{\partial y} \frac{\partial d}{\partial y} + \int_{\Omega}d $$\bar{u}$$ \cdot \nabla c  +\int_{\Omega} \mu c d  $\\ 

\hspace{0.7cm} l(d)=$\int_{\Omega} q d $\\
\begin{remark}
The bilinear form a(c,d) is not symmetric i.e. a(c,d) $\neq$ a(d,c).
\end{remark}

\subsection{Existence and uniqueness of the solution}
For a non-symmetric bilinear form $a(\cdot,\cdot)$, Lax Milgram theorem \cite{RefK} guarantees both existence and uniqueness of the solution of variational form .The theorem is stated as follows:\\
\begin{theorem}
The variational problem $a(c,d)=l(d) \hspace{0.5mm} \forall d \in V$ has unique solution $c \in V$ if the following conditions hold\\
(i)$a(\cdot,\cdot)$ is continuous i.e. $\exists$ $N_a$ such that $a(c,d) \leq N_a \|c\|_V \|d\|_V \forall c \in V , d \in V$\\
(ii)$a(\cdot,\cdot)$ is coercive i.e. $\exists$ $K_a$ such that $a(d,d) \geq K_a \|d\|_V^2 \forall d \in V$ and,\\
(iii) $l(\cdot)$ is continuous linear functional on V i.e. $\exists$ $N_l$ such that $l(d) \leq N_l \|d\|_V  \forall d \in V$ \\
\end{theorem}

Let us show that these conditions hold for this variational form (3) too, so that we can conclude the existence and uniqueness of its solution.\vspace{2mm}\\
(i) \textit{Continuity of $a(\cdot,\cdot)$}\\
It is clear from the above bilinear form,
\begin{equation}
    a(c,d) \leq \mid \int_{\Omega}D_x\frac{\partial c}{\partial x}\frac{\partial d}{\partial x} \mid + \mid \int_{\Omega}D_y\frac{\partial c}{\partial y} \frac{\partial d}{\partial y} \mid + \mid  \int_{\Omega}d \bar{u} \cdot \nabla c \mid + \mid \int_{\Omega} \mu c d \mid
\end{equation}
We will find the bounds for each term separately. Since the calculation for the first two terms will be same, therefore it is enough to calculate the bound for one of them. Before proceeding let us introduce few relations between norms: \vspace{1mm}\\
$ \|d\|$ $\leq$ $(\|d\|^2+\| \frac{\partial d}{\partial x} \|^2 + \| \frac{\partial d}{\partial y} \|^2)^\frac{1}{2}$ = $\|d\|_1 , \forall d \in V$\\
Similarly $ \|\frac{\partial d}{\partial x}\|$ $\leq$  $\|d\|_1 $ and $ \|\frac{\partial d}{\partial y}\|$ $\leq$ $\|d\|_1 , \forall d \in V$ \vspace{2mm}\\
We start with the first term,
\begin{equation}
    \begin{split}
      \mid \int_{\Omega}D_x\frac{\partial c}{\partial x}\frac{\partial d}{\partial x} \mid & \leq \int_{\Omega} \mid D_x\frac{\partial c}{\partial x}\frac{\partial d}{\partial x} \mid \\
      & \leq D_1 \| \frac{\partial c}{\partial x}\| \| \frac{\partial d}{\partial x} \| \hspace{2mm} (H\ddot{o}lders \hspace{1mm}inequality)\\
      & \leq D_1 \| c\|_1 \|  d \|_1
    \end{split}
\end{equation}
where $D_1$ = $\underset{(x,y) \in \Omega}{sup} \mid  D_x \mid $ \vspace{2mm}\\
Similarly, \\
\begin{equation}
    \mid \int_{\Omega}D_y\frac{\partial c}{\partial y} \frac{\partial d}{\partial y} \mid  \leq D_2 \| c\|_1 \|  d \|_1
\end{equation}
where $D_2$ = $\underset{(x,y) \in \Omega}{sup} \mid D_y \mid$ \vspace{2mm}\\ 
Following the above steps the second and third term will be bounded as follow: 
\begin{equation}
     \begin{split}
        \mid  \int_{\Omega}d \bar{u} \cdot \nabla c \mid &= \mid \int_{\Omega}  d u \frac{\partial c}{\partial x}  +  \int_{\Omega}  d v \frac{\partial c}{\partial y} \mid \\
        & \leq \int_{\Omega} \mid d u \frac{\partial c}{\partial x} \mid +  \int_{\Omega} \mid d v \frac{\partial c}{\partial y} \mid \\
        & \leq (U+V) \|c\|_1 \|d\|_1
     \end{split}
 \end{equation}
 where U= $\underset{(x,y) \in \Omega}{sup} \mid u\mid$  , V=$\underset{(x,y) \in \Omega}{sup} \mid v\mid$ 
 \begin{equation}
     \begin{split}
        \mid \int_{\Omega} \mu c d \mid & \leq  \int_{\Omega}  \mid \mu c d \mid\\
        & \leq \mid \mu \mid \|c\|_1 \|d\|_1
     \end{split}
 \end{equation}
 Now putting these bounds into (4) 
 \begin{equation}
     a(c,d) \leq N_a \|c\|_1 \|d\|_1 \hspace{2mm}\forall c \in V, d \in V
 \end{equation}
 where $N_a = (D_1+D_2+U+V+\mid \mu \mid)$ \vspace{3mm}\\
 (ii) \textit{Coercivity of $a(\cdot,\cdot)$}\\
 \begin{equation}
     \begin{split}
         a(d,d) &= \int_{\Omega}D_x\frac{\partial d}{\partial x}\frac{\partial d}{\partial x}+ \int_{\Omega}D_y\frac{\partial d}{\partial y} \frac{\partial d}{\partial y} + \int_{\Omega}d \bar{u} \cdot \nabla d +\int_{\Omega} \mu  d^2\\
         &= \int_{\Omega}D_x(\frac{\partial d}{\partial x})^2 + \int_{\Omega}D_y(\frac{\partial d}{\partial y})^2  +  +\int_{\Omega} \mu  d^2\\
         & \geq D \|d\|_1^2
     \end{split}
 \end{equation}
 where the constant D is the minimum taken over all the coefficients on $\Omega$. The third term in the first equation vanishes by using incompressibility condition and Dirichlet boundary condition  on integrating. Denoting D by $K_a$, we have 
 \begin{equation}
    a(d,d) \geq K_a \|d\|_1^2 \hspace{2mm}\forall d \in V
 \end{equation}
 (iii) \textit{Continuity of $l(\cdot)$}
 \begin{equation}
    \begin{split}
        l(d) & \leq \hspace{1mm}  \mid \int_{\Omega} q d \mid \\
        & \leq \int_{\Omega} \mid q d \mid \\
        & \leq \|q\| \|d\| \\
        & \leq N_l \|d\|_1
    \end{split}
 \end{equation}
 where $N_l = \|q\|$\\
 Hence there exists unique solution $c \in V$ for this variational problem.
 
\subsection{Galerkin finite element formulation}
Let  $V_h$ $\subset$ V  be a suitable finite dimensional subspace of V. Let us introduce the finite element discretization of the domain:
we consider that the bounded domain $\Omega$ is discretized into element subdomains $\Omega_k$, with boundary $\partial \Omega_k$, for k=1,2,...,$n_{el}$, where $n_{el}$ is the number of elements. Let $h_k$ be the diameter of each element $\Omega_k$ and h =$\underset{k=1,2,...,n_{el}}{max}$ $h_k$ \vspace{2mm}\\
Now the Galerkin formulation is : Find $c_h \in V_h$ such that\\
\begin{equation}
    a(c_h,d_h) = l(d_h) \quad \forall d_h \in V_h
\end{equation}
Let $n_{pt}$ be the total number of nodes occurred after finite element discretization and $N^a$ be the standard basis function of node a, for a= 1,2,...,$n_{pt}$. Therefore the functions $c_h$ and $d_h$ belonging to $V_h$ can be interpolated as,\vspace{2mm}\\
$c_h$= $\sum_{a=1}^{n_{pt}} c_a N^a$ \quad $d_h$ = $\sum_{a=1}^{n_{pt}} d_a N^a $\\
Putting them into above finite element problem, we have\\
\begin{equation}
    \begin{split}
        a(\sum_{a=1}^{n_{pt}} c_a N^a,\sum_{b=1}^{n_{pt}} d_b N^b) & = l(\sum_{b=1}^{n_{pt}} d_b N^b)\\
        \sum_{a,b=1}^{n_{pt}} d_b a(N^a,N^b) c_a & = \sum_{b=1}^{n_{pt}} d_b l(N^b) \\
        D^t A C & = D^t L\\
    \end{split}
\end{equation}
where A is a matrix of order $n_{pt} \times n_{pt}$ with components $a(N^a,N^b)$, for a,b =1,2,...,$n_{pt}$ and L is a column vector with components $l(N^b)$, for b=1,2,...,$n_{pt}$\\
Hence the original discrete problem is equivalently expressed in the following linear system
\begin{equation}
    AC=L
\end{equation}
In this paper only the Galerkin method is considered.
\begin{remark}
Among various finite elements, we have considered here P1 elements whose shape functions are piecewise linear functions.
\end{remark}

\section{Stability estimates: Strong stability condition}
Strong stability condition for the exact solution plays an important role in finding a priori and a posteriori error estimates in the $L_2$ norm. Therefore before going to deduce these estimates it is important to establish the condition for this problem. \textit{This condition implies that the exact solution can be bounded by the source term}.\\
Let us start with taking $L_2(\Omega)$ norm over q and using equation (2),
\begin{equation}
    \begin{split}
        \|q\|^2 &= \int_{\Omega} (-\nabla \cdot \widetilde{\nabla}c + \overline{u} \cdot \nabla c + \mu c)^2\\
        &= \int_{\Omega} (-\nabla \cdot \widetilde{\nabla}c + \overline{u} \cdot \nabla c )^2 + \int_{\Omega} \mu^2 c^2 -2\mu \int_{\Omega} (\nabla \cdot \widetilde{\nabla}c)c + 2\mu \int_{\Omega}c\overline{u} \cdot \nabla c\\
        &=  \int_{\Omega} (-\nabla \cdot \widetilde{\nabla}c + \overline{u} \cdot \nabla c )^2 + \int_{\Omega} \mu^2 c^2 + 2\mu \int_{\Omega} (D_x (\frac{\partial c}{\partial x})^2 + D_y (\frac{\partial c}{\partial y})^2) 
    \end{split}
\end{equation}
The fourth term in the second line vanishes by using incompressibility condition along with Dirichlet boundary condition on integrating and the third term in the last line comes after integrating the third term in the previous line along with using the same boundary condition.\\
From the last equation we can conclude that\\
\begin{equation}
\begin{split}
     \|q\| & \geq \mu \|c\| \\
     \|q\| &\geq \sqrt{2\mu D_1'} \|\frac{\partial c}{\partial x}\|\\
     \|q\| &\geq \sqrt{2\mu D_2'} \|\frac{\partial c}{\partial y}\|
\end{split}
   \end{equation}
 where the positive constants $D_1', D_2'$ come from taking minimum values of the coefficients $D_x , D_y$ respectively over $\Omega$.\\
Now considering new constants $E_1, E_2, E_3$ just by taking inverse of the constants $\mu, \sqrt{2\mu D_1'}, \sqrt{2\mu D_2'}$ appearing in the equation (17) respectively, we can rewrite the above inequalities in the following form,
\begin{equation}
\begin{split}
     \|c\| & \leq E_1 \|q\| \\
     \|\frac{\partial c}{\partial x}\| &\leq E_2 \|q\|\\
     \|\frac{\partial c}{\partial y}\| &\leq E_3 \|q\|
\end{split}
   \end{equation}
   Since we are aiming to bound the $H^2$ norm of the exact solution by $L_2$ norm of source term, therefore our remaining work is to bound the $L_2$ norm of double derivative terms of the exact solution by $\|q\|$.\vspace{3mm}\\
Using triangle inequality,
\begin{equation}
    \begin{split}
        \|\nabla \cdot \widetilde{\nabla}c\| &\leq \|\overline{u} \cdot \nabla c\| + \| \mu c\| + \|q\|\\
        & \leq \|u \frac{\partial c}{\partial x}\| + \|v\frac{\partial c}{\partial y}\| + \| \mu c\| + \|q\|\\
        & \leq U  \|\frac{\partial c}{\partial x}\|+ V  \|\frac{\partial c}{\partial y}\| + \mid \mu \mid \|c\| + \|q\| \\
        & \leq E_4 \|q\|
    \end{split}
\end{equation}
where  the constant $E_4= 1+ E_1 \mid \mu \mid + E_2 U + E_3 V$; the constants U,V have the same meaning as in section (2.4).\vspace{2mm}\\
Again, expanding the terms in $\nabla \cdot \widetilde{\nabla}c$,
\begin{equation}
    \begin{split}
\|\nabla \cdot \widetilde{\nabla}c\|^2 &= \int_{\Omega} (D_x \frac{\partial^2 c}{\partial x^2}+ D_y \frac{\partial^2 c}{\partial y^2} + \frac{\partial D_x}{\partial x} \frac{\partial c}{\partial x} + \frac{\partial D_y}{\partial y} \frac{\partial c}{\partial y})^2\\
& \geq D'^2 \int_{\Omega}  (\frac{\partial^2 c}{\partial x^2}+  \frac{\partial^2 c}{\partial y^2} +  \frac{\partial c}{\partial x} +  \frac{\partial c}{\partial y})^2\\
& \geq \int_{\Omega}  (\frac{\partial^2 c}{\partial x^2}+  \frac{\partial^2 c}{\partial y^2})( \frac{\partial c}{\partial x} +  \frac{\partial c}{\partial y})\\
& = \int_{\partial\Omega}(\nabla c \cdot \bar{n})(\frac{\partial c}{\partial x} +  \frac{\partial c}{\partial y}) - \int_{\Omega} \nabla c \cdot \nabla (\frac{\partial c}{\partial x} +  \frac{\partial c}{\partial y})
    \end{split}
\end{equation}
where the positive constant $D'$ comes from taking minimum of all the coefficients over $\Omega$ and $\bar{n}$ is the outward normal vector to the boundary. Here we ignore the first term by considering it either positive or zero on the boundary for specific c.\\
Now, expanding the remaining term of the last line and further using H$\ddot{o}$lders' inequality for each term, we have,
\begin{equation}
    \begin{split}
    \int_{\Omega} \nabla c \cdot \nabla (\frac{\partial c}{\partial x} +  \frac{\partial c}{\partial y}) & = \int_{\Omega} (\frac{\partial c}{\partial x} \frac{\partial }{\partial x}(\frac{\partial c}{\partial x}+\frac{\partial c}{\partial y})+\frac{\partial c}{\partial y}\frac{\partial }{\partial x}(\frac{\partial c}{\partial x}+\frac{\partial c}{\partial y}))\\
    & = \int_{\Omega} (\frac{\partial c}{\partial x} \frac{\partial^2 c}{\partial x^2}  + \frac{\partial c}{\partial x} \frac{\partial^2 c}{\partial x \partial y} + \frac{\partial c}{\partial y} \frac{\partial^2 c}{\partial y \partial x} + \frac{\partial c}{\partial y}\frac{\partial^2 c}{\partial y^2})\\
     & \leq \| \frac{\partial c}{\partial x}\| \|\frac{\partial^2 c}{\partial x^2} \| + (\|\frac{\partial c}{\partial x}\|+\|\frac{\partial c}{\partial y}\|)  \|\frac{\partial^2 c}{\partial x \partial y} \| + \|\frac{\partial c}{\partial y}\| \|\frac{\partial^2 c}{\partial y^2} \| \\
    & \leq E_5 \|q\| \hspace{0.5mm} (\|\frac{\partial^2 c}{\partial x^2} \| + \|\frac{\partial^2 c}{\partial x \partial y} \| + \|\frac{\partial^2 c}{\partial y^2} \|)\\
    \end{split}
\end{equation}
where $E_5= 2 (E_2+E_3)$\\
Putting the above inequality in (20) and squaring both sides we will have,
\begin{equation}
\begin{split}
 \|\nabla \cdot \widetilde{\nabla}c\|^4  & \geq  E_5^2 \|q\|^2 \hspace{0.5mm} (\|\frac{\partial^2 c}{\partial x^2} \| + \|\frac{\partial^2 c}{\partial x \partial y} \| + \|\frac{\partial^2 c}{\partial y^2} \|)^2\\ 
 & \geq  E_5^2 \|q\|^2 (\|\frac{\partial^2 c}{\partial x^2} \|^2 + \|\frac{\partial^2 c}{\partial x \partial y} \|^2 + \|\frac{\partial^2 c}{\partial y^2} \|^2)
\end{split}
\end{equation}
Using the above result into equation (19) and letting $E_6^2= E_4^4/E_5^2$,
\begin{equation}
     \|\frac{\partial^2 c}{\partial x^2} \|^2 + \|\frac{\partial^2 c}{\partial x \partial y} \|^2 + \|\frac{\partial^2 c}{\partial y^2} \|^2 \leq E_6^2 \|q\|^2
\end{equation}
By definition of $H^2(\Omega)$ norm, we  have,
\begin{equation}
    \begin{split}
        \|c\|_2^2 &= \|c\|^2 +\| \frac{\partial c}{\partial x} \|^2 + \|\frac{\partial c}{\partial y} \|^2 + \|\frac{\partial^2 c}{\partial x^2} \|^2 + \|\frac{\partial^2 c}{\partial x \partial y} \|^2 + \|\frac{\partial^2 c}{\partial y^2} \|^2\\
        & \leq C_s^2 \|q\|^2
    \end{split}
\end{equation}
where $C_s^2= E_1^2+ E_2^2+E_3^2+E_6^2$.\\
Since norm is a positive quantity, hence the strong stability estimate in its final form is:
\begin{equation}
    \boxed{\|c\|_2\leq C_s \|q\|}
\end{equation}
\section{A priori estimate in the energy norm}
We can define energy norm if the operator is positive definite. This property holds for this operator  as $\mu$ is positive. Now the energy norm is defined by:\\

$\|d\|_E^2 = a(d,d) \hspace{2mm} \forall d \in V$ \vspace{2mm}\\
In this non-symmetric operator the diffusion term introduces symmetric part of the operator where as the advection term represents skew symmetric part. Therefore the bilinear form can be systematically decomposed into these two parts, viz.\\

a(c,d) = $a_{symm}(c,d)$ + $a_{skew}(c,d)$ \vspace{2mm}\\
where

$a_{symm}(c,d)$ = $\frac{1}{2}$(a(c,d) $+$ a(d,c))= $\int_{\Omega}D_x\frac{\partial c}{\partial x}\frac{\partial d}{\partial x}+ \int_{\Omega}D_y\frac{\partial c}{\partial y} \frac{\partial d}{\partial y} +\int_{\Omega} \mu c d $\\

$a_{skew}(c,d)$ = $\frac{1}{2}$ (a(c,d) $-$ a(d,c)) =  $\int_{\Omega}d $$\bar{u}$$ \cdot \nabla c$ \vspace{2mm}\\
The skew symmetric part does not contribute to the energy norm as, $a_{skew}(d,d)$ = 0 ; hence\\

 $\|d\|_E^2$ = $a_{symm}(d,d)$ \vspace{2mm}\\
 For the \textbf{symmetric} part of the bilinear form, the Cauchy-Schwarz inequality holds with respect to the above norm and is given by\\
\begin{equation}
    \mid  a_{symm} (c,d) \mid \leq \|c\|_E \|d\|_E
\end{equation}
The dual norm approach \cite{RefB} has been considered to accomodate the skew symmetric part. The canonical definition of the 'skew norm' is given by, \\

$\|c\|_{skew}$ = $\underset{d \in V}{sup}$ $\frac{\mid a_{skew}(c,d) \mid}{\|d\|_E}$ \vspace{2mm}\\
Clearly the above definition implies, 
\begin{equation}
    \mid a_{skew}(c,d) \mid \leq  \|c\|_{skew} \|d\|_E
\end{equation}
Let e be the finite element error, $\tilde{c}_h$ be the nodal interpolant of the exact solution c and $c_h$ be the finite element solution.\vspace{3mm}\\
e = $c_h$ - c \vspace{1mm}\\
  = $c_h$ - $\tilde{c}_h$ + $\tilde{c}_h$ - c \vspace{1mm}\\
   = $e_h$ + $\eta_c$ \vspace{3mm}\\
  where, $e_h$ = ($c_h$ - $\tilde{c}_h$) is the portion of the error in $V_h$ and $\eta_c$ = ($\tilde{c}_h$ - c) is the interpolation errror belonging to V. \vspace{2mm}\\
 Using the concept of Galerkin orthogonality we have, $a(e,e_h)$ = 0 for $e_h$ $\in$ $V_h$ and considering $\mid a_{skew}(c,d) \mid$ = $\mid a_{skew}(d,c) \mid$\\
 \begin{equation}
  \begin{split}
   \|e\|_E^2&= \mid a(e,e) \mid\\
            &= \mid a(e,e_h + \eta_c) \mid \quad by \hspace{1mm} error \hspace{1mm} splitting\\ 
            &= \mid a(e,\eta_c) \mid \quad as \hspace{1mm} a(e,e_h) = 0\\
            &= \mid a_{symm}(e,\eta_c) + a_{skew}(e,\eta_c) \mid\\
            & \leq \mid a_{symm}(\eta_c,e) \mid + 
            \mid a_{skew}(\eta_c,e) \mid  \\
            & \leq \|\eta_c\|_E \|e\|_E + \|\eta_c\|_{skew} \|e\|_E \quad using \hspace{1mm}(26)\hspace{1mm} and\hspace{1mm} (27)
 \end{split}   
 \end{equation}
 Therefore,
 \begin{equation}
  \|e\|_E  \leq \|\eta_c\|_E + \|\eta_c\|_{skew} 
  \end{equation}
 Standard interpolation estimate \cite{RefB} is of the form,
 \begin{equation}
     \|\eta_c\|_1 \leq \bar{C}(p,\Omega) h^p \|c\|_{p+1}
 \end{equation}
The exact solution c has been assumed of regularity r$\geq$ p+1.\\
 Now,
 \begin{equation}
\begin{split}
\|\eta_c\|_E^2& = \mid a_{symm}(\eta_c,\eta_c) \mid\\
            & = \mid \int_{\Omega}D_x\frac{\partial \eta_c}{\partial x}\frac{\partial \eta_c}{\partial x}+ \int_{\Omega}D_y\frac{\partial \eta_c}{\partial y} \frac{\partial \eta_c}{\partial y} +\int_{\Omega} \mu \eta_c^2 \mid\\
            & \leq C_1'^2\mid \int_{\Omega}(\frac{\partial \eta_c}{\partial x})^2+ \int_{\Omega}(\frac{\partial \eta_c}{\partial y})^2  +\int_{\Omega} \eta_c^2 \mid\\
            & = C_1'^2 \| \eta_c \|^2_1
     \end{split}
 \end{equation}
 where the constant $C_1'$ comes after taking maximum of all the coefficients over $\Omega$. Using (30) and letting $\bar{C}_1$=$C_1'$ $\bar{C}(p,\Omega)$,
 \begin{equation}
     \|\eta_c\|_E \leq \bar{C}_1 h^p \|c\|_{p+1}
 \end{equation}
 Now,
 \begin{equation}
 \begin{split}
\|\eta_c\|_{skew} &= \underset{d \in V}{sup} \frac{\mid B_{skew}(\eta_c ,d) \mid}{\|d\|_E}\\
                 &= \underset{d \in V}{sup} \frac{\mid \int_{\Omega} d \bar{u} \cdot \nabla \eta_c \mid }{(\int_{\Omega}D_x(\frac{\partial d}{\partial x})^2+ \int_{\Omega}D_y(\frac{\partial d}{\partial y})^2 +\int_{\Omega} \mu d^2)^\frac{1}{2}}
 \end{split}
 \end{equation}
 Now we will find bounds for numerator and denominator separately.
\begin{equation}
    \begin{split}
      \int_{\Omega} d \bar{u} \cdot \nabla \eta_c \mid & \leq  \int_{\Omega} \mid d u \frac{\partial \eta_c}{\partial x} \mid +  \int_{\Omega} \mid d v \frac{\partial \eta_c}{\partial y} \mid \\
     & \leq U \|d\| \| \frac{\partial \eta_c}{\partial x}\| + V \|d\| \|\frac{\partial \eta_c}{\partial y}\| \quad using \hspace{1mm} H\ddot{o}lder's \hspace{1mm} inequality\\
     & \leq (U+V) \|d\|_1 \|\eta_c\|_1\\
    \end{split}
\end{equation}
where U,V and the norm inequalities are same as mentioned in section (2.4) and the denominator,
\begin{equation}
    \begin{split}
  (\int_{\Omega}D_x(\frac{\partial d}{\partial x})^2+ \int_{\Omega}D_y(\frac{\partial d}{\partial y})^2 +\int_{\Omega} \mu d^2)^\frac{1}{2} & \geq D' (\int_{\Omega}(\frac{\partial d}{\partial x})^2+ \int_{\Omega}(\frac{\partial d}{\partial y})^2  +\int_{\Omega} d^2)^\frac{1}{2}\\
  & = D' \|d\|_1 
\end{split}
\end{equation}
where the constant $D'=D^{1/2}$ which is the same constant appeared in section (2.4) in the proof of coercivity of the bilinear form.\\
Substituing the results obtained from (34) and (35) into (33) and letting $C_2'$= $\frac{U+V}{D}$, we will have
\begin{equation}
 \| \eta_c\|_{skew}  \leq  C_2' \hspace{1mm} \| \eta_c\|_1
\end{equation}
Using (30) and letting $\bar{C}_2$=$C_2'$ $\bar{C}(p,\Omega)$, (36) becomes
\begin{equation}
     \|\eta_c\|_{skew} \leq \bar{C}_2 h^p \|c\|_{p+1}
 \end{equation}
 Substituing the results obtained from (32) and (37) into (29) and letting $C'$=$C_1' + C_2'$, our a priori error estimate in energy norm will be
 \begin{equation}
 \boxed{  \|e\|_E \leq C' h^p \|c\|_{p+1}}
 \end{equation}
 provided the exact solution c has regularity r $\geq$ p+1.
 \begin{remark}
 For piecewise linear elements i.e. for P1 elements the above inequality will be
 \begin{equation}
     \|e\|_E \leq C' h \|c\|_2
 \end{equation}
 Now by applying strong stability condition,
 \begin{equation}
     \|e\|_E \leq C' C_s h \|q\|
 \end{equation}
 \end{remark}
 
 \section{A posteriori error estimate in the energy norm}
In this method the bilinear form does not need to be split into symmetric and non-symmetric parts. The first few lines of the derivation are same as above. Here we have deduced residual based a posteriori error estimation which is expressed in terms of residuals and hence computable.\\
Using error splitting and Galerkin orthogonality, introduced in the previous section we get,
\begin{equation}
    \begin{split}
        \|e\|_E^2 &= \mid a(e,\eta_c) \mid\\
                  &= \mid a(c_h-c, \eta_c) \mid \\
                  &= \mid a(c_h, \eta_c) - a(c, \eta_c) \mid \\
                  &= \mid a(c_h, \eta_c) - l(\eta_c) \mid 
\end{split}
\end{equation}
We have obtained the last line by using the weak form (3), and thus the expression gets rid of the exact solution. Now the job, is to form residuals, which involves integration by parts. This has to be done separately on each element $\Omega_k$, as the finite element solution $c_h$ is only continuous across the element interfaces.\\
Introducing the notation\\

$\int_{\Omega} := \sum_{k=1}^{n_{el}} \int_{\Omega_k}$ , $\int_{\partial \Omega'} := \sum_{k=1}^{n_{el}} \int_{\partial \Omega_k}$ \vspace{2mm}\\
where $\Omega$ is the given domain.
\begin{equation}
    \begin{split}
        \mid a(c_h, \eta_c) - l(\eta_c) \mid &= \mid \int_{\Omega}D_x\frac{\partial c_h}{\partial x}\frac{\partial \eta_c}{\partial x}+ \int_{\Omega}D_y\frac{\partial c_h}{\partial x} \frac{\partial \eta_c}{\partial y} + \int_{\Omega} \eta_c \bar{u} \cdot \nabla c_h +\int_{\Omega} \mu c_h \eta_c - \int_{\Omega} q \eta_c \mid\\
        &= \mid \sum_{k=1}^{n_{el}} (\int_{\Omega_k} \widetilde{\nabla}c_h \cdot \nabla \eta_c + \int_{\Omega_k} \eta_c \bar{u} \cdot \nabla c_h + \int_{\Omega_k} \mu c_h \eta_c - \int_{\Omega_k} q \eta_c) \mid \\
        &= \mid \sum_{k=1}^{n_{el}} (\int_{\partial \Omega_k} (\widetilde{\nabla}c_h \cdot \bar{n})\eta_c - \int_{\Omega_k} (\nabla \cdot \widetilde{\nabla} c_h) \eta_c + \int_{\Omega_k} (\eta_c \bar{u} \cdot \nabla c_h +  \mu c_h \eta_c -  q \eta_c) )\mid\\
        &= \mid \int_{\partial \Omega'} (\widetilde{\nabla}c_h \cdot \bar{n})\eta_c + \int_{\Omega'}r_h \eta_c \mid \\
        & = \mid \sum_{k=1}^{n_{el}} \int_{\Omega_k} r_h \eta_c \mid \\ 
        & \leq \sum_{k=1}^{n_{el}} \int_{\Omega_k} \mid r_h \eta_c \mid
        \end{split}
\end{equation}
where $\bar{n}$ denotes unit outward normal. The first term of the fourth line vanishes assuming $\eta_c$ zero on inter element boundaries. \\
The residual, $r_h$ = $-\bigtriangledown \cdot \widetilde{\bigtriangledown}c_h + \overline{u} \cdot \bigtriangledown c_h + \mu c_h - q $\\
Before going further, let us introduce appropriate interpolation estimates, given in \cite{RefA} viz.\\
\begin{equation}
    \|\eta_c\|_{L_2 (\Omega_k)} \leq C_{I,k} h_k \|e\|_{E_,\Omega_k} 
\end{equation}
where
\begin{equation}
    C_{I,k}= \underset{d \in V}{sup}  \frac{h_k^{-1} \|\tilde{d_h}-d\|_{L_2 (\Omega_k)}}{ \|d\|_{E_,\Omega_k}}
\end{equation}
where $\tilde{d_h}$ is the interpolant of d. By standard interpolation theory the constant $C_{I,k}$ is bounded. \vspace{2mm}\\
Now we will find the bound for element interior terms obtained in (42) using H$\ddot{o} $lders inequality and the appropriate interpolation estimate, mentioned above, as follows:
\begin{equation}
   \begin{split}
    \sum_{k=1}^{n_{el}} \int_{\Omega_k} \mid r_h \eta_c \mid  &\leq \sum_{k=1}^{n_{el}} \| r_h\|_{L_2 (\Omega_k)} \|\eta_c\|_{L_2 (\Omega_k)}\\
    & \leq \sum_{k=1}^{n_{el}} \| r_h\|_{L_2 (\Omega_k)} C_{I,k} h_k \|e\|_{E_,\Omega_k}\\
    & \leq (\sum_{k=1}^{n_{el}} C_{I,k} h_k \| r_h\|_{L_2 (\Omega_k)} ) \|e\|_E \\
    \end{split}
\end{equation}
Substituing the above results into (42) and equating it with (41), we will have the following form after knocking out the common positive term $\|e\|_E$ from both sides,
\begin{equation}
    \boxed{\|e\|_E \leq \sum_{k=1}^{n_{el}} C_{I,k} h_k \| r_h\|_{L_2 (\Omega_k)} }
\end{equation}
\begin{remark}
Clearly this a posteriori error estimate includes the computed solution instead of exact solution, hence it is computable with appropriate estimates of $C_{I,k}$, for k = 1,2,...,$n_{el}$.
\end{remark}

\section{A priori error estimation in the $L_2$ norm}
Here we are going to use Nitsche trick to find first a priori $L_2$ estimate and later a posteriori error estimate in $L_2$ norm too.\\
\subsection{Nitsche trick}
In general Nitsche trick involves an auxiliary problem, here it is dual problem, following \cite{RefB}. The dual problem consists of the adjoint operator with the error as the source term. Now the variational form of the dual problem including homogeneous Dirichlet boundary conditions is stated as follows: \\
Find $\beta \in $ V such that $\forall$ d $\in$ V
\begin{equation}
    a(d,\beta)=(d,e)
\end{equation}
Replacing d by e in the above equation
\begin{equation}
    (e,e)=\|e\|^2 = a(e,\beta)
\end{equation}
Applying the strong stability condition, obtained in (25) to the adjoint problem, we will have
\begin{equation}
    \|\beta\|_2 \leq C_s \|e\|
\end{equation}
Let $\beta_h \in V_h$ be the interpolant of $\beta$. By the Galerkin orthogonality 
\begin{equation}
    a(e,\beta_h) = 0
\end{equation}
Subtracting (50) from (48) and denoting the interpolation error, $(\beta - \beta_h)$ by $\eta_\beta$, we have
\begin{equation}
    \begin{split}
        \|e\|^2 &= a(e,\beta-\beta_h)\\
                &= \mid a(e, \eta_\beta) \mid\\
                & \leq \mid a_{symm}(e,\eta_\beta) \mid + \mid a_{skew}(e,\eta_\beta) \mid \\
                &= \mid a_{symm}(\eta_\beta , e) \mid + \mid a_{skew}(\eta_\beta, e) \mid \\
                & \leq  \|\eta_\beta\|_E \|e\|_E + \|\eta_\beta\|_{skew} \|e\|_E \\
                &\leq (\| \eta_\beta\|_E + \| \eta_\beta\|_{skew}) C' h^p \|c\|_{p+1} 
    \end{split}
\end{equation}
The last line is obtained by applying a priori error estimation in the energy norm, obtained in (38). \vspace{2mm}\\
The interpolation estimate \cite{RefB} for $\eta_\beta$ is:
\begin{equation}
    \| \eta_\beta\|_1  \leq \bar{C}(p,\Omega) h \| \beta\|_2 \\
\end{equation}
Replacing $\| \beta\|_2$  by applying strong stability condition (49) to the above inequality, we have
\begin{equation}
     \| \eta_\beta\|_1 \leq \bar{C_s} h \|e\|
\end{equation}
where $\bar{C_s}=C_s \bar{C}(p,\Omega)$ \\
We can obtain the bounds for $\|\eta_\beta\|_E$ and $\|\eta_\beta\|_{skew}$ in the similar way done in the section (4), just by substituting $\eta_\beta$ in place of $\eta$. Using those bounds obtained in (31) and (36) in the equation (51), we will have
\begin{equation}
    \begin{split}
        \|e\|^2 & \leq ( (C_1'+C_2') \|\eta_\beta\|_1 ) C' h^p \|c\|_{p+1}\\
                & \leq (C_1'+C_2') \bar{C_s} C' h^{p+1} \|e\| \|c\|_{p+1} \\
        \|e\| & \leq C'' h^{p+1} \|c\|_{p+1}        
    \end{split}
\end{equation}
where $C''$= $(C_1'+C_2') C_s C'$
\begin{remark}
Nitsche trick extracts the extra power of h, therefore whereas the convergence in the energy norm is $O(h^p)$, here it is of $O(h^{p+1})$.
\end{remark}
\begin{remark}
In particular for P1 elements the estimation will be
\begin{equation}
    \|e\| \leq C'' C_s h^2 \|q\| \hspace{3mm} (by \hspace{1mm} stability \hspace{1mm} estimate)
\end{equation}
\end{remark}

\section{A posteriori estimates in the $L_2$ norm}
We will start with Nitsche trick as done in the previous section. Further we will follow another way to get rid of the exact solution. In the procedure we will use some results deduced in section (5).

\subsection{Nitsche trick}
Let us again introduce the variational form in the similar way as follows:\\
Find $\beta \in $ V such that $\forall$ d $\in$ V
\begin{equation}
    a(d,\beta)=(d,e)
\end{equation}
where the error e is the source term.
Recall the strong stability condition mentioned by the equation (49), as it has very important role in obtaining error estimate in $L_2$ norm.\\
Substituing e in the place of d into the above equation,
\begin{equation}
    \begin{split}
        \|e\|^2 &= (e,e)\\
              &= a(e, \beta)\\
              &= a(c_h-c, \beta)\\
              &= a(c_h, \beta) - a(c, \beta)\\
              &= a(c_h, \beta) - l(\beta)
    \end{split}
\end{equation}
where $c_h$ is the finite element solution and last line comes from the weak form (3).\\
Let $\beta_h \in V_h$ be the interpolant of $\beta$. By Galerkin formulation (13), we will have
\begin{equation}
    a(c_h,\beta_h)-l(\beta_h) = 0
\end{equation}
Subtracting (58) from (57) and following the same notation $\eta_\beta$ for the interpolation error,  $(\beta-\beta_h)$, we arrive at the following form: 
\begin{equation}
\begin{split}
    \|e\|^2 & = a(c_h,\beta-\beta_h)-l(\beta-\beta_h)\\
    & \leq \mid a(c_h, \eta_\beta)-l(\eta_\beta) \mid\\
    & \leq  \sum_{k=1}^{n_{el}} \int_{\Omega_k} \mid r_h \eta_\beta \mid \\
\end{split}
    \end{equation}
    The last inequality comes from the derivation, done in section (5) to obtain the result (42), just by substituting $\eta$ by $\eta_\beta$, where $r_h$ is the residual term.\\
The appropriate interpolation estimates for each element \cite{RefA} are as follows:
\begin{equation}
   \|\eta_\beta\|_{L_2 (\Omega_k)} \leq C_{i,k} h_k^2\|\beta\|_2  
\end{equation}
where  the constants $C_{i,k}$, for k=1,2,...,$n_{el}$ depend upon the their respective subdomains.\\
Let us find the bound for element interior terms using the strong stability condition (49) and the above interpolation estimates as follows:
\begin{equation}
   \begin{split}
    \sum_{k=1}^{n_{el}} \int_{\Omega_k} \mid r_h \eta_\beta \mid  &\leq \sum_{k=1}^{n_{el}} \| r_h\|_{L_2 (\Omega_k)} \|\eta_\beta\|_{L_2 (\Omega_k)}\\
    & \leq \sum_{k=1}^{n_{el}} \| r_h\|_{L_2 (\Omega_k)}C_{i,k} h_k^2\|\beta\|_2\\
    & \leq (\sum_{k=1}^{n_{el}} C_s C_{i,k} h_k^2 \| r_h\|_{L_2 (\Omega_k)} ) \|e\| \\
    \end{split}
\end{equation}
Substituting the result obtained in (61) into (59) and cancelling out the common positive term $\|e\|$ from both side, we will remain with
\begin{equation}
    \boxed{\|e\| \leq C_s \sum_{k=1}^{n_{el}}  C_{i,k} h_k^2 \| r_h\|_{L_2 (\Omega_k)} }
\end{equation}

\section{Numerical experiment}
In this section we assess the second order numerical convergence of the finite element scheme advocated by error estimate. To test the theoretically established order of convergence numerically we consider the following two test cases from literature in Hydrology \cite{RefF}:\textit{first case} is only for homogeneous Dirichlet boundary condition and \textit{second case} is for mixed boundary condition. In both the cases exact solution is known.
\subsection{First case}
We have taken a square bounded domain (0,1)x(0,1) as $\Omega$ on which the incompressible fluid flow considered. This flow respects ADR equation with spatially variable coefficients along with homogeneous Dirichlet boundary condition. We have taken the values of the coefficients from \cite{RefF}, which is motivated by the fact of real-life application from hydrology. We have taken care of the fact that we are dealing with an incompressible flow. We have used \textit{freefem++} to solve the problem and to find error in $L_2$ norm and its order of convergence. Let us first introduce the values of the coefficients we have worked with.\\
u= 0.5 (1+0.02x),v= -0.5(1+0.02y),$D_x= 0.2 (1+0.02x)^2$, $D_y=0.02(1+0.02y)^2$ and $\mu$= 0.01\\
The exact solution $c=sin(x(x-1)y(y-1))$ on $\Omega$ along with c= 0 on $\partial \Omega$\vspace{1mm}\\
\begin{remark}
Table 1. shows the error and order of convergence obtained in $L_2$ norm for this homogeneous Dirichlet boundary condition.
\end{remark}



\subsection{Second case} In this case the horizontal parallel boundaries are treated as non-homogeneous Dirichlet boundaries with c(y=0)=1 and c(y=1)=1. On the rest of the boundaries homogeneous Neumann condition have been considered. Velocity components and diffusion coefficients are set at values similar to those in first case. \\
The exact solution c = cos(x(x-1)y(y-1)) on $\Omega$ \vspace{1mm}\\
\begin{remark}
Table 2 shows the error and order of convergence obtained in $L_2$ norm for mixed boundary conditions, which is mentioned in this case.
\end{remark}

 
 \begin{table}
   \centering
    \begin{tabular}{||c c c||}
    \hline 
    Mesh size & Error in $L_2$ norm &  Order of convergence \\
    \hline \hline
      10      &  0.000242665        &                       \\
      \hline
      20      &  5.72783 $e^{-5}$   & 2.08291  \\
      \hline
      40      &  1.27668 $e^{-5}$   & 2.16559  \\
      \hline
      80      &  3.15467 $e^{-6}$   & 2.01684  \\
      \hline
      160     &  8.22607 $e^{-7}$   & 1.93921  \\
      \hline
      320     &  1.81488 $e^{-7}$   & 2.18033 \\
      \hline
    \end{tabular}
     \caption{Error and Order of convergence obtained in $L_2$ norm for the first case}
    \label{table:draglift1}
\end{table}

\begin{table}[]
    \centering
    \begin{tabular}{||c c c||}
    \hline 
    Mesh size & Error in $L_2$ norm &  Order of convergence \\
    \hline \hline
      10      &  1.55274 $e^{-5}$        &                       \\
      \hline
      20      &  7.55422 $e^{-6}$   & 1.03946  \\
      \hline
      40      &  1.64313 $e^{-6}$   & 2.20083  \\
      \hline
      80      &  3.80341 $e^{-7}$   & 2.11108  \\
      \hline
      160     &  9.92829 $e^{-8}$   & 1.93768  \\
      \hline
      320     &  2.81933 $e^{-8}$   & 1.81619 \\
      \hline
    \end{tabular}
\caption{Error and Order of convergence obtained in $L_2$ norm for the second case}
    \label{table:draglift2}
\end{table}

\begin{remark}
It is clear from the tables appended above that in both the cases the order of convergence for the computed results is approximately 2, which coincides with our theoretical approximation.
\end{remark}
\section{Conclusion}
Both a priori and a posteriori error estimates of ADR equation with variable coefficients in $L_2$ norm lead to second order convergence of the finite element scheme and this has been numerically verified.


\begin{thebibliography}{}
\bibitem{RefA}
C.Johnson and P.Hansbo, \textit{Adaptive finite element methods in computational mechanics}, Comput. Methods Appl. Mech. Engrg., 101(1992), 143-181.
\bibitem{RefB}
J.R.Stewart and T.J.R. Hughes, \textit{A tutorial in elementary finite element error analysis: A systematic presentation of a priori and a posteriori error estimates}, Comput. Methods Appl. Mech. Engrg. 158(1998), 1-22
\bibitem{RefC}
R.Verf$\ddot{u}$rth, \textit{A posteriori error estimators for convection-diffusion equations}, Numer. Math., 80(1998), 641-663.
\bibitem{RefD}
C.E.Baumann and J.T.Oden, \textit{A discontinuous hp finite element method for convection-diffusion problems}, Comput. Methods Appl. Mech. Engrg., 175(1999), 311-341.
\bibitem{RefE}
M.Vohral$\Acute{i}$k, \textit{A posteriori error estimates for lower-order mixed finite element discretizations of convection-diffusion-reaction equations}, SIAM J. Numer. Anal.,Vol.45, No-4,(2007), 1570-1599.
\bibitem{RefF}
A.Sanskrityayn, V.P.Singh, V.K.Bharati and N.Kumar, \textit{Analytical solution of two-dimensional advection-dispersion equation with spatio-temporal coefficients for point sources in an infinite medium using Green's function method}, Environ. Fluid Mech., 18(2018), 739-757.
\bibitem{RefG}
L.Angermann and Erlangen, \textit{Balanced a posteriori error estimates for finite-volume type discretizations of convection-dominated elliptic problems}, Computing, 55(1995), 305-323.
\bibitem{RefH}
A.Ern, A.F.Stephansen and M.Vohral$\Acute{i}$k, \textit{Guaranteed and robust discontinuous Galerkin a posteriori error estimates for convection-diffusion-reaction problems}, J. Comput. Appl. Math., 234(2010), 114-130.
\bibitem{RefI}
P.Houston, C.Schwab and E.S$\ddot{u}$li, \textit{Discontinuous hp-finite element methods for advection-diffusion-reaction problems}, SIAM J. Numer. Anal., Vol.39, No.6 (2002),2133-2163.
\bibitem{RefJ}
R.Verf$\ddot{u}$rth, \textit{Robust a posteriori error estimators for a singularly perturbed reaction-diffusion equation}, Numer. Math., 78(1998), 479-493.
\bibitem{RefK}
D.Bahuguna, V.Raghvendra and B.V.R.Kumar, \textit{Topics in sobolev spaces and applications}, Narosa Publishing House, 2002, ISBN: 8173194297.


\end{thebibliography}
\end{document}